\documentclass[runningheads]{llncs}

\usepackage{graphicx}
\usepackage{hyperref}

\usepackage{amsmath}
\usepackage{amssymb}
\usepackage{bbm}
\usepackage{subcaption}
\usepackage[usenames, dvipsnames]{color}
\usepackage{wrapfig}
\usepackage{booktabs}
\usepackage{colortbl}

\definecolor{Gray}{gray}{0.95}
\renewcommand{\mathbb}{\mathbbm}

\newcommand{\scp}[2]{{\left\langle {#1}\, , \, {#2}\right\rangle}}
\newcommand{\lform}[2]{{\left( {#1}\, \left\vert \vphantom{#1}\,  {#2}\right.\right)}}

\renewcommand{\phi}{\varphi}
\renewcommand{\epsilon}{\varepsilon}

\graphicspath{{figures/}}

\begin{document}

\title{A Model for Elastic Evolution \\on Foliated Shapes}

\author
{
Dai-Ni Hsieh\inst{1} \and
Sylvain Arguill\`ere\inst{2} \and
Nicolas Charon\inst{1} \and \\
Michael I. Miller\inst{3} \and
Laurent Younes\inst{1}
}

\authorrunning{D.-N. Hsieh et al.}

\institute
{
Department of Applied Mathematics and Statistics, \\ Johns Hopkins University, Baltimore, MD 21218, USA \and
Institut Camille Jordan, Villeurbanne, France \and
Department of Biomedical Engineering, \\ Johns Hopkins University, Baltimore, MD 21218, USA
}

\maketitle


\begin{abstract}

We study a shape evolution framework in which the deformation of shapes from time $t$ to $t + dt$ is governed by a regularized anisotropic elasticity model. More precisely, we assume that at each time shapes are infinitesimally deformed from a stress-free state to an elastic equilibrium as a result of the application of a small force. The configuration of equilibrium then becomes the new resting state for subsequent evolution. The primary motivation of this work is the modeling of slow changes in biological shapes like atrophy, where a body force applied to the volume represents the location and impact of the disease. Our model uses an optimal control viewpoint with the time derivative of force interpreted as a control, deforming a shape gradually from its observed initial state to an observed final state. Furthermore, inspired by the layered organization of cortical volumes, we consider a special case of our model in which shapes can be decomposed into a family of layers (forming a ``foliation"). Preliminary experiments on synthetic layered shapes in two and three dimensions are presented to demonstrate the effect of elasticity.


\end{abstract}


\section{Introduction}

Understanding changes in anatomical shapes is an essential problem for the analysis of many diseases, and more specifically for their tracking over time. While a large variety of shape analysis methods have been successful in cross-sectional problems, in which one exhibits differences in some shape spaces between two populations, a longitudinal analysis of shape changes in a single subject requires a  careful modeling of the process, including the structure of the deforming material and the process causing the changes.

Our interest is more specifically in the modeling of dynamical changes in biological tissues that typically occur over long periods of time (such as the progressive atrophy of an organ along the years of disease) thus resulting in sequences of slow alterations in the material properties. Unlike several past approaches for longitudinal analysis based on generic deformation analysis pipelines \cite{QIU2009S51,Durrleman2013}, this work constitutes a step towards a model that could explain morphological variations through more physically interpretable features, while also enabling the estimation of the potential source location and severity of the pathology.    

The setting we introduce in this paper is in part inspired from elastic models describing slow changes in biological shapes \cite{rodriguez1994stress,lubarda2002mechanics,dicarlo2002growth,amar2005growth,tallinen2016growth}. We propose a  shape evolution framework in which the transition between time $t$ and $t+dt$ is governed by a regularized anisotropic elasticity model in which shapes at rest (or stress free) at time $t$ find a new elastic equilibrium at time $t+dt$ as a result of the application of a small force $dF$, reaching a configuration that becomes their new resting state for subsequent evolution.

The paper is organized as follows. Sect.~\ref{sec:general} introduces a general mathematical setting for the shape evolution process, while Sect.~\ref{sec:elastic} details the specific laminar elastic model used in our experiments. Sect.~\ref{sec:experiments} provides some numerical illustrations of the method, focusing, so far, on synthetic data. We conclude the paper in Sect.~\ref{sec:discussion}.


\section{Shape Evolution Paradigm}
\label{sec:general}


We first introduce a general shape evolution model which we shall make specific to layered elastic materials in the next section. Given an initial shape $M_0$ which we take as a compact domain of $\mathbb{R}^3$, we consider deformed shapes $t \mapsto \varphi(t, M_0)$ in which $\varphi$ is a time-dependent diffeomorphism obtained as the flow of a time-dependent velocity field $t \mapsto v(t, M_0)$  satisfying the system:
\begin{align}
	\left\{
		\begin{aligned}
			& \partial_t \hspace{0.5pt} \varphi(t) = v(t, \varphi(t)), \ \varphi(0) = \mathrm{id} \\
			& L_{\varphi(t)} v(t) = j(t) \, dp
		\end{aligned}
	\right. .
	\label{eq:syst.1}
\end{align}
We shall assume that each $v(t)$ belongs to $V$, a reproducing kernel Hilbert space (RKHS) of vector fields on $\mathbb{R}^3$, continuously embedded in the space $C^2_0(\mathbb R^3, \mathbb R^3)$ which denotes the Banach space of $C^2$ vector fields on $\mathbb{R}^3$ such that $u$, $Du$, and $D^2 u$ vanish at infinity, equipped with the norm defined by $\|u\|_{2, \infty} = \|u\|_\infty + \|Du\|_\infty + \|D^2 u\|_\infty$. ($Du$ denotes the differential of $u$.) Also, in~\eqref{eq:syst.1}, $L_\varphi$ is a certain linear operator from $V$ to $V^*$, $j(t)$ is a vector field on $\mathbbm R^3$ supported by $\varphi(t, M_0)$, and $j(t) \, dp$ is its associated element in $V^*$. In the rest of the paper, for any generic function $(t,x) \mapsto f(t,x)$, we adopt the convention of writing $f(t)$ instead of $f(t, \cdot)$.

Our objective consists in estimating $j(t)$ given an observed initial shape $M_0$ and target shape $M_1$, under the evolution model expressed in \eqref{eq:syst.1}. In this general setting, such an inverse problem is typically formulated as an optimal control problem by adding either regularization penalties or constraints on $j$ for well-posedness.  

In spite of the fact that $\phi(t)$ is a diffeomorphism obtained as the flow a vector field, the important difference between this approach and standard registration of volumetric data in diffeomorphic frameworks such as \cite{beg2005computing} is that the evolution here is not controlled directly by the velocity $v$. Instead, $v$ is implicitly determined by $j$ through the operator $L_\varphi$. The definition of $L_{\phi}$ thus considerably affects the nature and properties of optimal solutions. In the next section, we define $L_\varphi$ for layered structures based on anisotropic linear elasticity.

\section{Elastic Evolution of Layered Shapes}
\label{sec:elastic}
We build  the operator $L_\varphi$ step by step. First, for a fixed elastic shape with elastic tensor $\Lambda$, we define the operator $\mathcal{L}_\Lambda$ such that the solution of $\mathcal{L}_\Lambda u = F \, dp$ characterizes, up to a small regularization, the elastic displacement from resting state to equilibrium when a force density $F$ is applied. Next, we define an elastic tensor $\Lambda$ in the special case of layered shapes. We then specify how the elastic tensor is transformed to $\Lambda^\varphi$ when a diffeomorphism $\varphi$ acts on layered shapes, which is needed to keep track of the configuration as it gets progressively deformed.

By our assumption that the deformed shapes are at rest at each time, the infinitesimal displacement at each time is given by $\mathcal{L}_{\Lambda^\varphi} \delta u = \delta F \, dp$ for some infinitesimal force density $\delta F$. After normalization by $\delta t$, we finally obtain the elastic operator for layered shapes as $L_\varphi = \mathcal{L}_{\Lambda^\varphi}$ such that $\mathcal{L}_{\Lambda^\varphi} v = j \, dp$.

\subsection{Linear Elastic Model for Small Deformation}

Let $M \subset \mathbb{R}^3$ be a compact domain, which is assumed to represent an elastic material at rest. When applying a small force $F$ to $M$, the displacement field $u$ of the small deformation $\psi = id + u$ can be obtained based on the minimum total potential energy principle, in which $u$ minimizes
\[
U_\Lambda(u) - \int_M F(p)^\top u(p) \, dp\,,
\]
where $U_\Lambda(u)$ is the linear elastic energy resulting from the displacement $u$ on $M$, and $\Lambda$ is the elastic tensor describing elastic properties of the material. The exact expression of $U_\Lambda(u)$ is detailed below. Here, we assume that only a body force ``inside the volume'', i.e., a force density, affects the material, but that no 
pressure acts on its boundary.

In order to ensure, eventually, a diffeomorphic evolution of the shape when passing to a time-dependent model, we add a regularization term to the total potential energy and look for a regularized response $u$ to the force $F$ given by
\begin{align}
	u
	=
	\underset{u' \in V}{\arg\min} \left(
	\frac{\delta}{2} \, \|u'\|_V^2 + U_\Lambda(u') - \int_M F(p)^\top u'(p) \, dp \right) ,
	\label{eq:u_min}
\end{align}
in which $\delta >0$ is a small regularization parameter.

Following \cite{gurtin1973linear,marsden1994mathematical,ciarlet1988three}, we  now specify the definition of the elastic energy and characterize the solution of~\eqref{eq:u_min}. For a displacement field $u$ on $M$, we denote the linear strain tensor by
\[
	\varepsilon(u) = \frac{1}{2} \, (Du + Du^\top) . 
\]
Let $\varepsilon_p(u)$ denote the evaluation of $\varepsilon(u)$ at a point $p \in M$. The linear elastic energy corresponding to $u$ can then be written as
\[
	U_\Lambda(u) = \int_M \Lambda_p(\varepsilon_p(u))\, dp ,
\]
where the elastic tensor $p \mapsto \Lambda_p$ is a mapping from $M$ to the set of non-negative quadratic forms over the space of $3$ by $3$ symmetric matrices. The function $\Lambda$ encodes the elastic properties of the material. For example, an isotropic linear elastic material has $\Lambda_p(\varepsilon) = \frac{\lambda}{2} \left( \sum_{i = 1}^3 \varepsilon_{ii} \right)^2 + \mu \sum_{i, j = 1}^3 \varepsilon_{ij}^2$ for all $p \in M$, where $\lambda$ and $\mu$ are the Lam{\'e} parameters. 

If we further assume that $p\mapsto \|\Lambda_p\|$ is integrable over $M$, then there exists a linear operator $\mathcal A_\Lambda: C^1_0(\mathbb R^3, \mathbb R^3)\to C^1_0(\mathbb R^3, \mathbb R^3)^*$ uniquely defined by 
\[
	U_\Lambda(u) = \frac12 \lform{\mathcal A_\Lambda u}{u} ,
\]
where $\lform{\alpha}{u}$ denotes the application of a linear form $\alpha$ to a vector $u$. Since $V \hookrightarrow C^2_0(\mathbb{R}^3, \mathbb{R}^3)$, $\mathcal A_\Lambda$ can also be seen as an operator from $V$ to $V^*$. Letting $L: V\to V^*$ denote the duality operator on $V$ such that $\scp{v}{w}_V = \lform{Lv}{w}$, one can define the operator $\mathcal L_\Lambda: V \to V^*$ by $\mathcal L_\Lambda =  \delta L + \mathcal A_\Lambda$ and then write
\[
\frac12 \lform{\mathcal L_\Lambda u}{u} = \frac\delta2 \, \|u\|_V^2 + U_\Lambda(u)\,.
\]
As a result, the solution of \eqref{eq:u_min} is given by $u$ such that $\mathcal L_\Lambda u = F \, dp$.

\subsection{Layered Shapes and Tangential Isotropic Elasticity}

\begin{figure}
	\centering
	\begin{subfigure}{0.3\textwidth}
		\centering
		\includegraphics[height = 80pt]{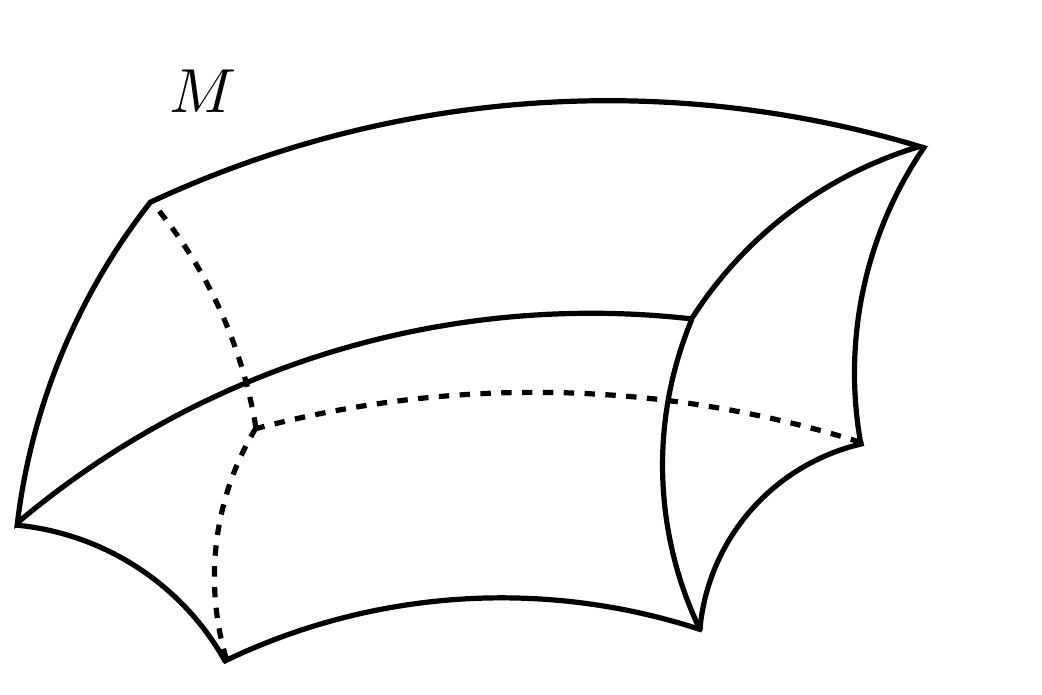}
		\caption{3D volume}
	\end{subfigure}
	\begin{subfigure}{0.4\textwidth}
		\centering
		\includegraphics[height = 80pt]{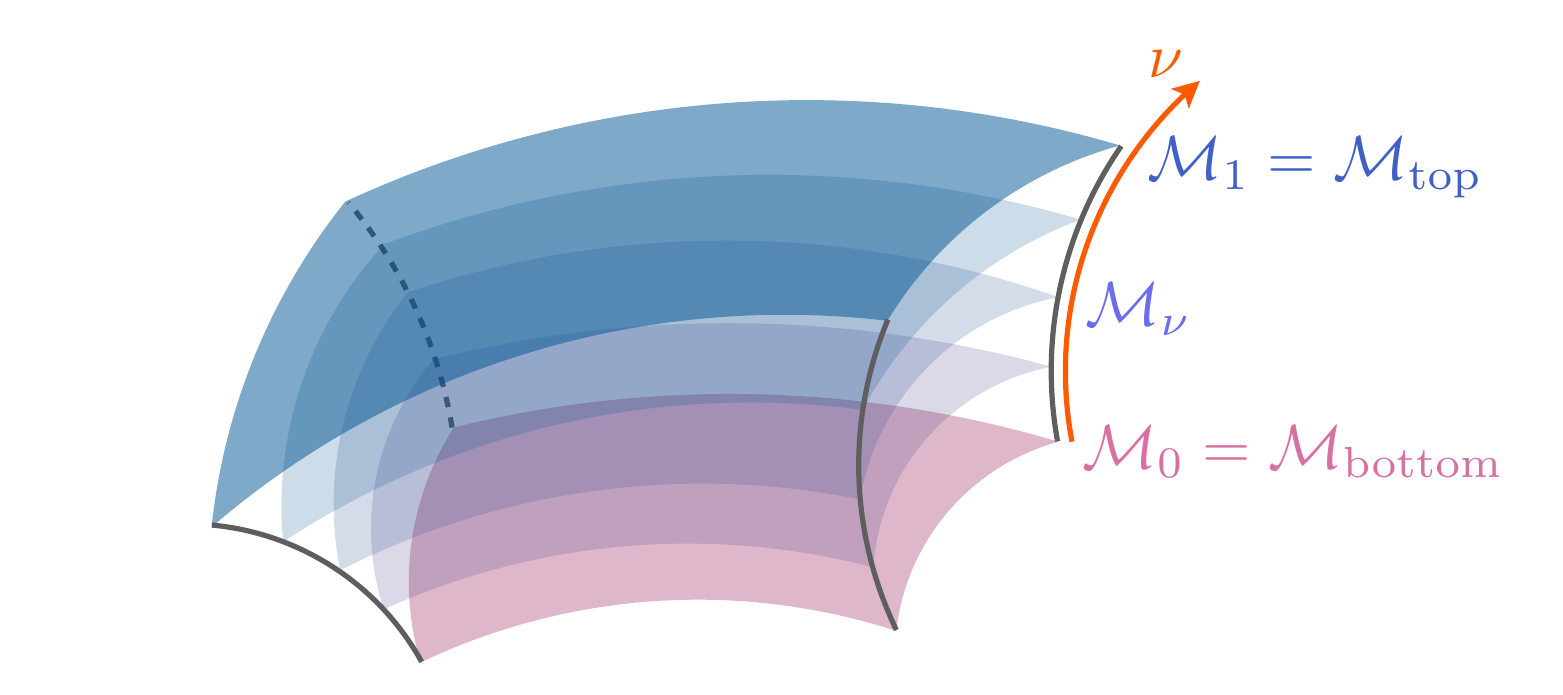}
		\caption{Foliation}
	\end{subfigure}
	\caption{Shapes are assumed to have a layered foliation.}
	\label{fig:foliation}
\end{figure}

Besides being elastic, we assume that $M$ comes with a layered structure as illustrated in Fig.~\ref{fig:foliation} and that the elastic behavior of $M$ is ``isotropic tangentially to layers.'' This is defined precisely below by specifying the particular form of $\Lambda$ in this configuration. One motivating example for such a construction is that of the cerebral cortex and its organization along cortical layers and columns.

For the compact domain $M \subset \mathbb{R}^3$, we assume that there are two surfaces $\mathcal M_{\mathrm{bottom}}$, the bottom layer, $\mathcal M_{\mathrm{top}}$, the top layer, included in $\partial M$. Moreover, we assume that the top layer is obtained from a diffeomorphism $\Phi: [0, 1] \times \mathcal M_{\mathrm{bottom}} \rightarrow M$, such that $\Phi(0) = \mathrm{id}$ and $\Phi(1,\mathcal M_{\mathrm{bottom}}) = \mathcal M_{\mathrm{top}}$. We let $S$ denote the vector field defined by $\partial_\nu \Phi(\nu, p) = S(\Phi(\nu, p))$, and we assume that $S$ is continuously differentiable.

Define intermediate layers as $\mathcal M_\nu = \Phi(\nu, \mathcal M_{\mathrm{bottom}})$, so that $\{\mathcal M_\nu\}_{\nu \in [0, 1]}$ forms a foliation of $M$ to which $S$ is {\em transversal}. Algorithms building such foliations of volumes have been introduced in brain imaging in order to estimate thickness, see e.g., \cite{jones2000three,das2009registration,ratnanather20183d}. Given such a transversal vector field $S$ and layers $\{\mathcal M_\nu\}_{\nu \in [0, 1]}$, we now define an elastic tensor $\Lambda$ which is consistent with this layered structure and ``isotropic on layers". For $p\in M$, let $\{T_1(p), T_2(p)\}$ be {any} orthonormal basis of $T_p \mathcal M_\nu$ (the tangent plane to $\mathcal M_\nu$ at $p$). We define the skew linear change of coordinates at $p$ as the $3$ by $3$ matrix $C_p = \left[T_1(p) \ T_2(p) \ \frac{S(p)}{|S(p)|} \right]$. Now we consider $\Lambda$ of the form $\Lambda_p(\varepsilon) = \bar\Lambda(C_p^\top \varepsilon \, C_p)$, where
\begin{align}
\bar \Lambda(\xi)
	&=
	\frac{1}{2} \, \mu_{\mathrm{tan}} \left( \xi_{11} + \xi_{22} \right)^2
	+
	\frac{1}{2} \, \lambda_{\mathrm{tan}} \left( \xi_{11}^2 + \xi_{22}^2 + 2 \, \xi_{12}^2 \right)
	+
	\frac{1}{2} \, \lambda_{\mathrm{tsv}} \, \xi_{33}^2
    \notag
	\\
	&\hspace{10pt}
	+
	\lambda_{\mathrm{ang}} \left( \xi_{13}^2 + \xi_{23}^2 \right) .
	\label{eq:3D_layered_strain_energy}
\end{align}
$\bar \Lambda$ is independent of the position $p$ if the parameters $\mu_\mathrm{tan}$, $\lambda_{\mathrm{tan}}$, $\lambda_{\mathrm{tsv}}$, and $\lambda_{\mathrm{ang}}$ are fixed, which is assumed in the following. Letting $\xi = C_p^\top \varepsilon_p \, C_p$, we note that $\xi_{11}, \xi_{22}$ measure the stretches along tangential directions $T_1(p), T_2(p)$ respectively, and $\xi_{33}$ measures the stretch along the transversal direction, which explains the particular form we take for $\bar \Lambda$.

Let $\{T_1'(p), T_2'(p)\}$ be another orthonormal basis of $T_p \mathcal M_\nu$. Then we have
\[
	C_p'
	=
	\left[
		\begin{array}{ccc}
			\\[-5pt]
			T_1'(p) & T_2'(p) & \displaystyle \frac{S(p)}{|S(p)|} \\[-5pt]
			\phantom{a}
		\end{array}
	\right]
	=
	\left[
		\begin{array}{ccc}
			\\[-5pt]
			T_1(p) & T_2(p) & \displaystyle \frac{S(p)}{|S(p)|} \\[-5pt]
			\phantom{a}
		\end{array}
	\right]
	\left[
		\begin{array}{cc}
			\\[-5pt]
			\,G(p)\  & 0 \, \\[5pt]
			\,0\     & 1 \,
		\end{array}
	\right]  ,
\]
where $G(p)$ is a $2$ by $2$ orthogonal matrix. It follows that $\mathrm{tr}(\xi)$, $\mathrm{tr}(\xi^2)$, and $\xi_{11}^2 + \xi_{22}^2 + 2 \, \xi_{12}^2$ are invariant under this transformation, so $\Lambda_p$ does not depend on the choice of orthonormal basis of $T_p \mathcal M_\nu$. $\Lambda$ can thus be thought as ``isotropic on $T_p \mathcal M_\nu$'' for all $p \in M$. \\

{\noindent \it Remark.} One can make the same construction for a 2D domain $M$, with a foliation $\{\mathcal{M}_\nu\}_{\nu \in [0, 1]}$ made of curves, instead of surfaces. One can then describe $T_p \mathcal M_\nu$ by a single unit vector $T(p)$ and define $C_p = \left[T(p) \ \frac{S(p)}{|S(p)|} \right]$, a 2 by 2 matrix. In this case, the elastic tensor simplifies to
\begin{align*}
	\bar\Lambda(\xi)
	&=
	\frac{1}{2} \, \lambda_{\mathrm{tan}} \, \xi_{11}^2
	+
	\frac{1}{2} \, \lambda_{\mathrm{tsv}} \, \xi_{22}^2
	+
	\lambda_{\mathrm{ang}} \, \xi_{12}^2 .
\end{align*}


\subsection{Boundary Condition}
As the elastic energy is insensitive to the effect of rigid motions, one usually needs to specify certain boundary conditions in order to ensure that the displacement $u$ resulting from $F$ is well-defined. Here, we impose boundary conditions on $u$ such as, for example, $u = 0$ on some subset of $M$ (typically the bottom layer $\mathcal{M}_{\mathrm{bottom}}$). Going back to the case of cortical volumes, this is consistent with the typical assumption that grey matter atrophy is balanced mainly by expansion of the cerebrospinal fluid (CSF) region \cite{Khanal2016}, which in our case corresponds the top layer of $M$. The bottom layer represents the gray/white matter boundary and should remain relatively stable. 

In the context of reproducing kernel Hilbert spaces, this boundary condition corresponds to replacing $V$ by its closed subspace
\[
	V_0 = \{u \in V: u_{|_{\mathcal M_{\mathrm{bottom}}}} \equiv 0\} ,
\]
which contains smooth vector fields that do not move $\mathcal M_{\mathrm{bottom}}$. Under the constraint, the displacement response $u$ to a force $F$ is now given by
\[
	u = \underset{u' \in V_0}{\arg\min} \left( \frac{\delta}{2} \left\| u' \right\|_V^2 + \frac12 \lform{\mathcal{A}_\Lambda u'}{u'} - \lform{F dp}{u'}\right)\,. 
\]
Because $V_0$ is a closed subspace of $V$, the induced norm makes it into an RKHS. Note that $\mathcal{A}_\Lambda$ can then be restricted as an operator from $V_0$ to $V_0^*$, so the same formal analysis applies,  leading to an operator from $V_0$ to $V_0^*$ that we still denote by $\mathcal{L}_\Lambda$.


\subsection{Action of Diffeomorphisms}

If the volume $M$ is transformed by a diffeomorphism $\varphi$ to $\varphi(M)$, the deformed layered structure can be specified by
\[
	\widetilde \Phi: [0, 1] \times \varphi(\mathcal{M}_{\mathrm{bottom}}) \rightarrow \varphi(M), \ \ 
	(\nu, \widetilde p) \mapsto \varphi(\Phi(\nu, \varphi^{-1}(\widetilde p))) .
\]
In particular, the transversal vector field $S$ becomes $\widetilde S$ such that $\widetilde S(\phi(p)) = D\phi(p) \, S(p)$. Let $\widetilde{\mathcal{M}}_\nu = \varphi(\mathcal{M}_\nu)$ and $\{\widetilde{T}_1(\widetilde p), \widetilde{T}_2(\widetilde p)\}$ be an orthonormal basis of $T_{\widetilde p} \widetilde{\mathcal{M}}_\nu$. One can then define a skew linear transformation $\widetilde C = \left[\widetilde T_1, \widetilde T_2, \frac{\widetilde S}{|\widetilde S|}\right]$ and the transformed elastic tensor $\Lambda^\varphi$ given by $\Lambda^\varphi_{\widetilde p}(\varepsilon) = \bar \Lambda(\widetilde{C}^\top_{\widetilde p} \varepsilon \,\widetilde{C}_{\widetilde p})$ with $\bar \Lambda$ unchanged. This directly provides operators $\mathcal A_{\Lambda^\phi}$ and $\mathcal L_{\Lambda^\phi} = L_\varphi$.


\subsection{Evolution Model for Large Deformation}

We now have all required elements in place to operate the large deformation model in \eqref{eq:syst.1} for the particular operator $L_{\varphi}$ defined in the previous sections. The connection between the small deformation model and the evolution equations can be made by considering that, at a given time $t$, the displacement $\delta u$ resulting from an infinitesimal force density $\delta F$ is given by $L_\varphi \hspace{1pt} \delta u = \delta F \, dp$. The resulting velocity is  $v = \frac{\delta u}{\delta t}$, and letting $j = \frac{\delta F}{\delta t}$ yields $L_\varphi \, v = j \, dp$.

The regularization in equation \eqref{eq:u_min} ensures that $v(t)$ is smooth enough such that, under some conditions controlling the size of the control $j(t)$, it generates a flow of diffeomorphisms. This approach is inspired by methods in image and shape registration that use similar regularization, such as the large deformation diffeomorphic metric mapping algorithm \cite{beg2005computing,younes2010shapes}.


\section{Numerical Results}
\label{sec:experiments}

In this section, we  showcase several numerical experiments as applications of our model. In all these experiments, we consider a simplified case where $j$ is of the form $j(t,\varphi(t,x)) = D\varphi(t,x) \, j_0(x)$, that is, $j$ is completely specified by its value $j_0$ at time zero. We describe in the next section a low dimensional model for $j_0$, specified by a small number of parameters. We then consider the inverse problem of estimating these parameters based on the observation of the original and deformed shapes.

\subsection{Force Model}
As a simple model of atrophy, we assume that $j_0$ is the gradient of a curved Gaussian which is defined by
\[
	g(p; c, h, \sigma_{\mathrm{tan}}, \sigma_{\mathrm{tsv}})
	=
	h 
	\exp\!
	\left(
		-\frac{d_{\mathrm{tan}}^2(p, c)}{2 \sigma_{\mathrm{tan}}^2}
		-\frac{d_{\mathrm{tsv}}^2(p, c)}{2 \sigma_{\mathrm{tsv}}^2}
	\right) ,
\]
where
\begin{wrapfigure}[12]{r}{0.4\textwidth}
	\vspace{-35pt}
	\centering
	\includegraphics[width = 0.4\textwidth]{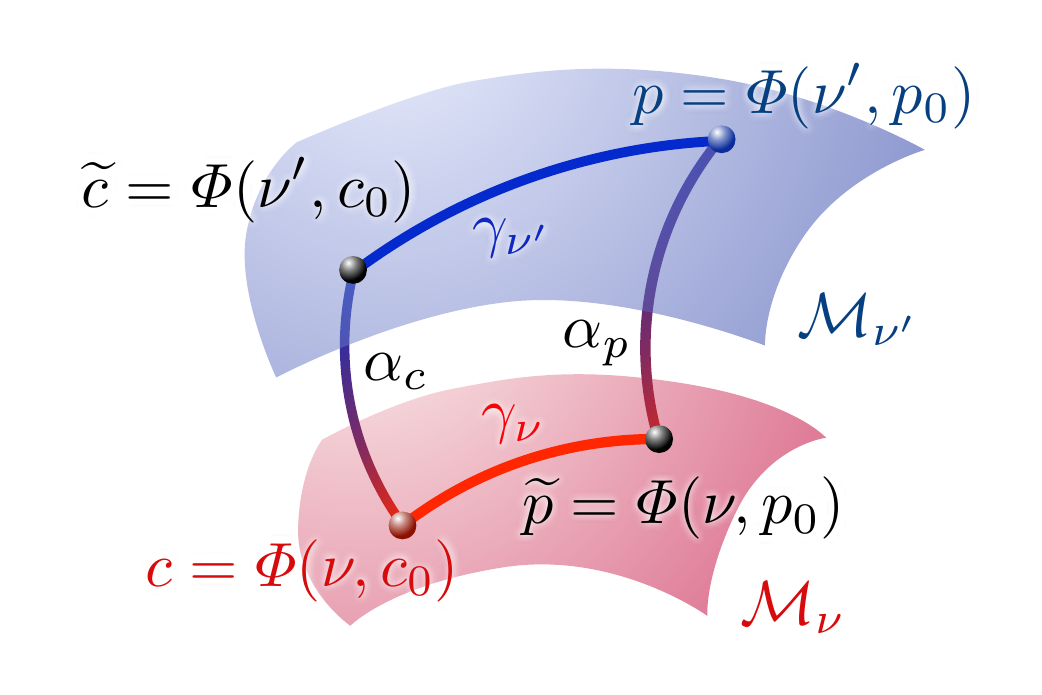}
	\caption{The four curves $\gamma_\nu$, $\gamma_{\nu'}$, $\alpha_c$, and $\alpha_p$ in the definition of curved Gaussian.}
	\label{fig:curvedGaussian_distances}
\end{wrapfigure}
\[
	d_{\mathrm{tan}}(p, c) = \frac{1}{2} \left( \mathrm{length}(\gamma_{\nu}) + \mathrm{length}(\gamma_{\nu'}) \right)
\]
and
\[
	d_{\mathrm{tsv}}(p, c) = \frac{1}{2} \left( \mathrm{length}(\alpha_c) + \mathrm{length}(\alpha_p) \right) .
\]
As shown in Fig.~\ref{fig:curvedGaussian_distances}, $\gamma_\nu$ and $\gamma_{\nu'}$ are geodesics on $\mathcal{M}_\nu$ and $\mathcal{M}_{\nu'}$ joining $c$, $\widetilde p$ and $\widetilde c$, $p$ respectively, and $\alpha_c$ and $\alpha_p$ are integral curves joining $c$, $\widetilde c$ and $\widetilde p$, $p$ respectively.

\subsection{Inverse Problem}

We assume that $\sigma_{\mathrm{tan}}$ and $\sigma_{\mathrm{tsv}}$ are known parameters and that the center of Gaussian $c$ is constrained on the middle layer $\mathcal{M}_{\frac{1}{2}}$. This assumption is partially motivated by the case of Alzheimer's disease in which histological observations suggest that tau protein accumulation preceding the onset occurs in internal cortical layers \cite{braak1991neuropathological}. Under all these assumptions, the control $j$, and, as a consequence, the entire evolution, is completely determined by the location of $c$ and the height of peak $h$. Given an undeformed shape $M_0$ and a target shape $M_1$, the inverse problem can now be stated as minimizing
\[
	J(c, h) = d(\varphi(1, M_0), M_1)
 \text{\quad subject to\quad }
 	\left\{
 		\begin{aligned}
 		& \partial_t \varphi(t) = v(t, \varphi(t)) \\
 		& \mathcal{L}_{\Lambda^{\varphi(t)}} v(t) = j(t) \, dp \\
 		& j(t) = (D\varphi(t) \, j_0)\circ \varphi(t)^{-1} \\
 		& j_0 = \nabla g \\
 		& c \in \mathcal{M}_{\frac{1}{2}}, \ h > 0
 \end{aligned}
 	\right. .
\]
We use in our experiments
\[
	d(M, M')	= \mathrm{volume}(M \triangle M') = \int_{\mathbb{R}^3} \left( \mathbbm{1}_{M} - \mathbbm{1}_{M'} \right)^2 dx ,
\]
which computes the volume of non-overlapping region of two shapes. Here, we assume that one does not have access to the layered structure for $M_1$. Therefore, we do not include any information of layers in the discrepancy measure $d$. 

Since the constraint $c \in \mathcal{M}_{\frac{1}{2}}$ burdens an optimization algorithm, we consider the coordinate representation of $J$ in $c$. In other words, let $(U, \psi)$ be a local chart of $c \in \mathcal{M}_{\frac{1}{2}}$. We define $\widehat J: \psi(U) \times \mathbb{R}_+ \rightarrow \mathbb{R}$ by $\widehat J(\widehat c, h) = J(\psi^{-1}(\widehat c), h)$. Instead of minimizing the objective function $J$, we minimize equivalently its coordinate representation $\widehat J$. By moving the domain of the objective function from the curved space $\mathcal{M}_{\frac{1}{2}} \times \mathbb{R}_+$ to the Euclidean space, we can then utilize a derivative-free optimization algorithm with only box constraints. \\

{\noindent \it Remark.} Following the terminology introduced in \cite{gris2015sub}, our model specifies a  ``deformation module,'' that we  call ``elastic module'' in the following. Such modules provide a deformation mechanism (here represented by $j_0$) that both drives the shape evolution and is advected by it (see \cite{gris2015sub} for more details). The free parameters for these modules are the control variables $(c,h)$ with additional geometric parameters $\theta_f = (\sigma_{\mathrm{tan}}, \sigma_{\mathrm{tsv}})$ for the force and $\theta_e = (\delta, \mu_{\mathrm{tan}}, \lambda_{\mathrm{tan}}, \lambda_{\mathrm{tsv}}, \lambda_{\mathrm{ang}})$ for the elastic properties of the volume. One can also relate our construction to that provided in \cite{grenander2007pattern}, in which our choice for $j_0$ provides what is called a ``sink'' in the referenced paper.

\subsection{Discretization}

Recall the layered structure $\mathcal M_\nu = \Phi(\nu, \mathcal M_{\mathrm{bottom}})$ where $\Phi$ is a diffeomorphism. We use a discrete set of layers (hence letting $\nu$ be an integer) with consistent triangulations. More precisely, we assume that vertices in $\mathcal M_\nu$ are $(p_1^\nu, \ldots, p_N^\nu)$ (with $N$ independent of $\nu$) and faces (which are triples or integer in $\{1, \ldots, N\}$) are also independent of $\nu$. We also define transverse edges between $p_k^\nu$ and $p_k^{\nu+1}$ and let $S(p_k^\nu) = p_k^{\nu+1} - p_k^\nu$ represent the discrete version of the transversal vector field (such a representation is provided, for example, by the algorithm introduced in \cite{ratnanather20183d}). This provides a decomposition of the volume $M$ into triangular prisms. Without adding vertices, we further decompose the prisms into a consistent tetrahedral mesh using a method similar to that introduced in \cite{porumbescu2005shell}. 

For the evaluation of the objective function, we use Dijkstra's algorithm to compute the length of geodesics needed to define the curved Gaussian map. Discretizing $j \, dp$ as a weighted sum of Dirac functions supported by vertices, one can show that the solution of $\mathcal{L}_{\Lambda^{\varphi}} v = j \, dp$ is uniquely specified by the values of $v$ at vertices. This results in a large linear system (including the constraint that  $v=0$ on the bottom layer), which is solved using conjugate gradient.


\subsection{Simulated Deformation}
\label{subsec:numerical_moduli}

\begin{figure}
	\centering
	\begin{subfigure}{0.49\textwidth}
		\centering
		\includegraphics[width = \textwidth]{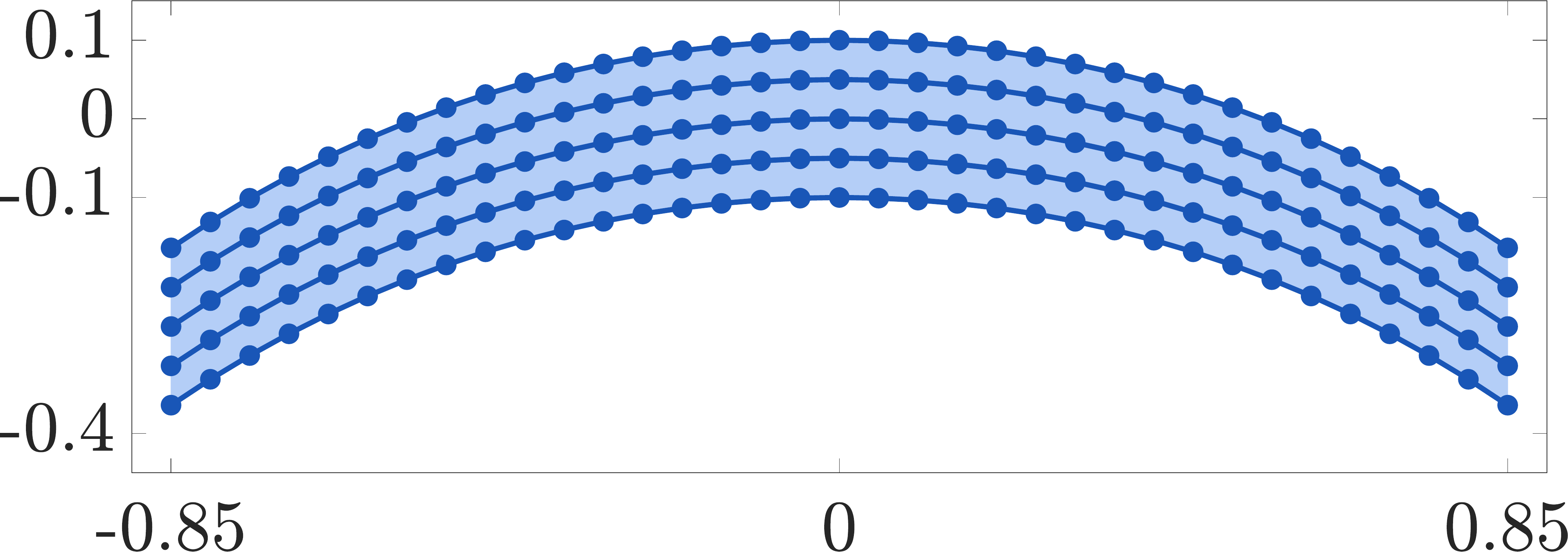}
		\caption{Undeformed shape}
		\label{fig:cap_template}
	\end{subfigure}
	\hfill
	\begin{subfigure}{0.49\textwidth}
		\centering
		\includegraphics[width = \textwidth]{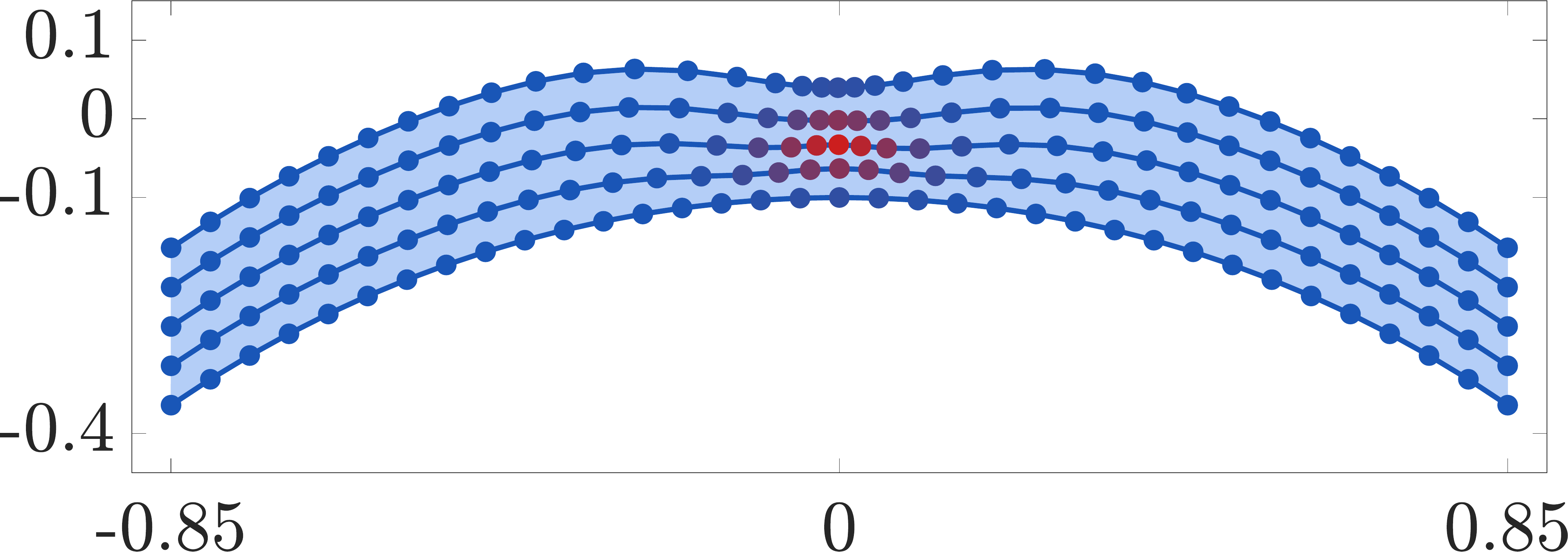}
		\caption{Tangent easy}
		\label{fig:cap_tan_easy}
	\end{subfigure}
	\\[10pt]
	
	\begin{subfigure}{0.49\textwidth}
		\centering
		\includegraphics[width = \textwidth]{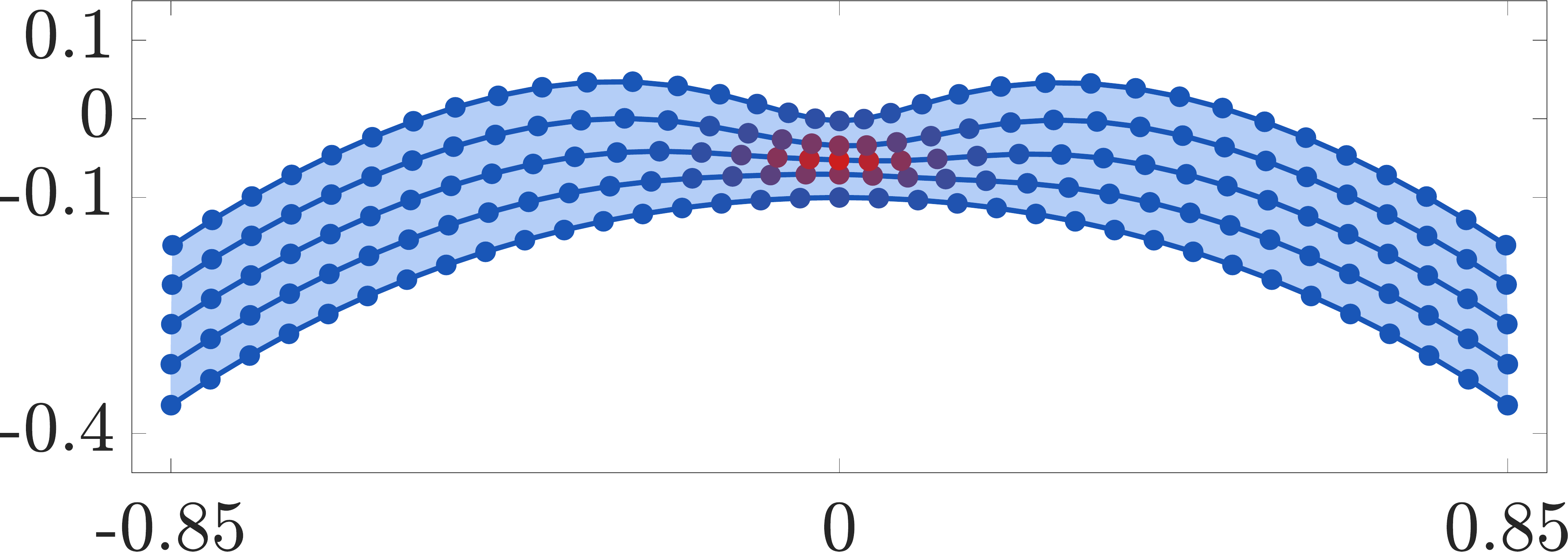}
		\caption{Transverse easy}
		\label{fig:cap_tsv_easy}
	\end{subfigure}
	\hfill
	\begin{subfigure}{0.49\textwidth}
		\centering
		\includegraphics[width = \textwidth]{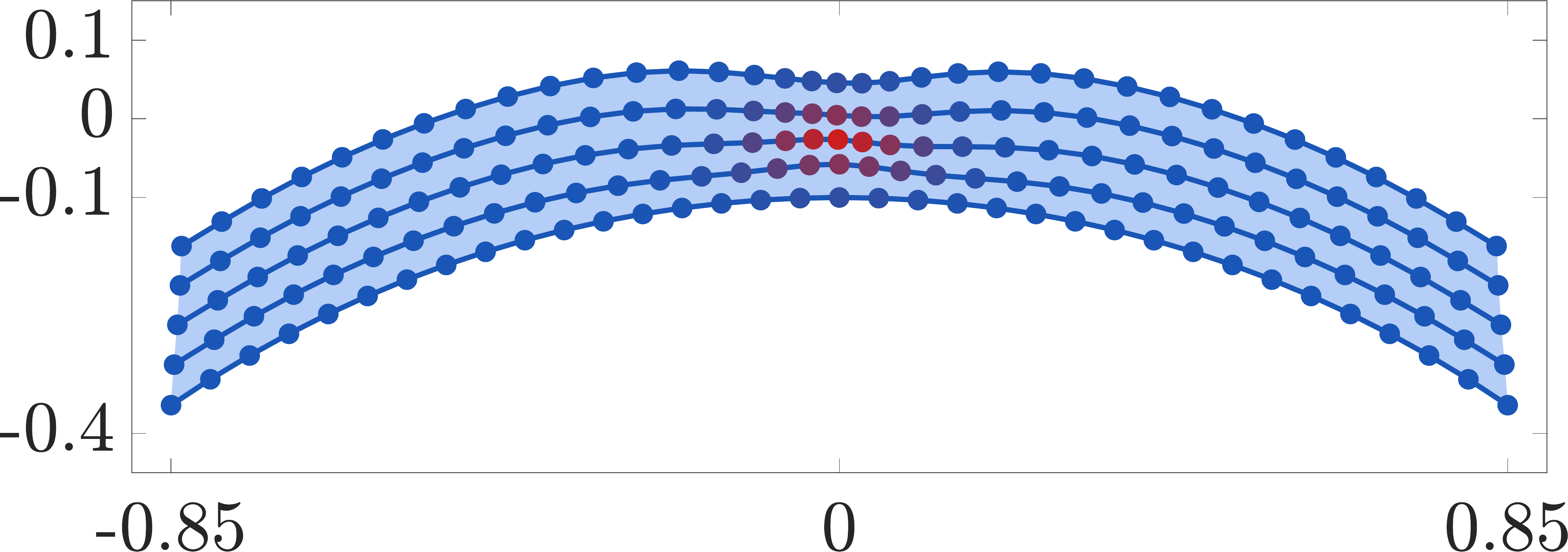}
		\caption{Angle easy}
		\label{fig:cap_ang_easy}
	\end{subfigure}
	\caption{Effects of different elastic parameters on shape deformations.}
\end{figure}

Figs.~\ref{fig:cap_tan_easy} to \ref{fig:cap_ang_easy} display the deformed shapes $\varphi(1, M_0)$ after applying the same $j_0$ to the undeformed shape $M_0$, which is shown in Fig.~\ref{fig:cap_template}. The colors of landmarks from blue to red represent the values of curved Gaussian from low to high. The force parameters are $c = (0, 0)$, $h = 3$, $(\sigma_{\mathrm{tan}}, \sigma_{\mathrm{tsv}}) = (0.1, 0.05)$. We used three different sets of elasticity parameters, that we call ``tangent easy,'' ``transverse easy'' and ``angle easy.'' All three cases use $V$ with a Mat\'ern kernel of order three and width 0.01, $\delta = 10^{-6}$,  $(\lambda_{\mathrm{tan}}, \lambda_{\mathrm{tsv}}, \lambda_{\mathrm{ang}}) = (1, 3, 3)$, $(3,1, 3)$ and $(3,3,1)$  in the tangent-easy, transverse-easy, and angle-easy cases, respectively.

The different elastic parameters produce different responses to the same $j_0$: in Fig.~\ref{fig:cap_tan_easy}, the landmarks mainly move tangentially; in Fig.~\ref{fig:cap_tsv_easy}, they mainly move transversely; while in Fig.~\ref{fig:cap_ang_easy}, we see the least movement due to the constraint on the change of lengths between landmarks. Recall that in all cases the bottom layer is kept fixed by our model. 

\subsection{Simulated Inverse Problem}
\label{subsec:numerical_2D}

We now present 2D and 3D simulated inverse problems, in which we first generated deformed targets using curved Gaussian, and then we tried to retrieve the parameters $c$ and $h$. In our experiments, we used \texttt{surrogateopt} in Matlab to find a minimizer.

\begin{figure}
	\centering
	\begin{subfigure}{0.43\textwidth}
		\begin{subfigure}{\textwidth}
			\centering
			\includegraphics[width = \textwidth]{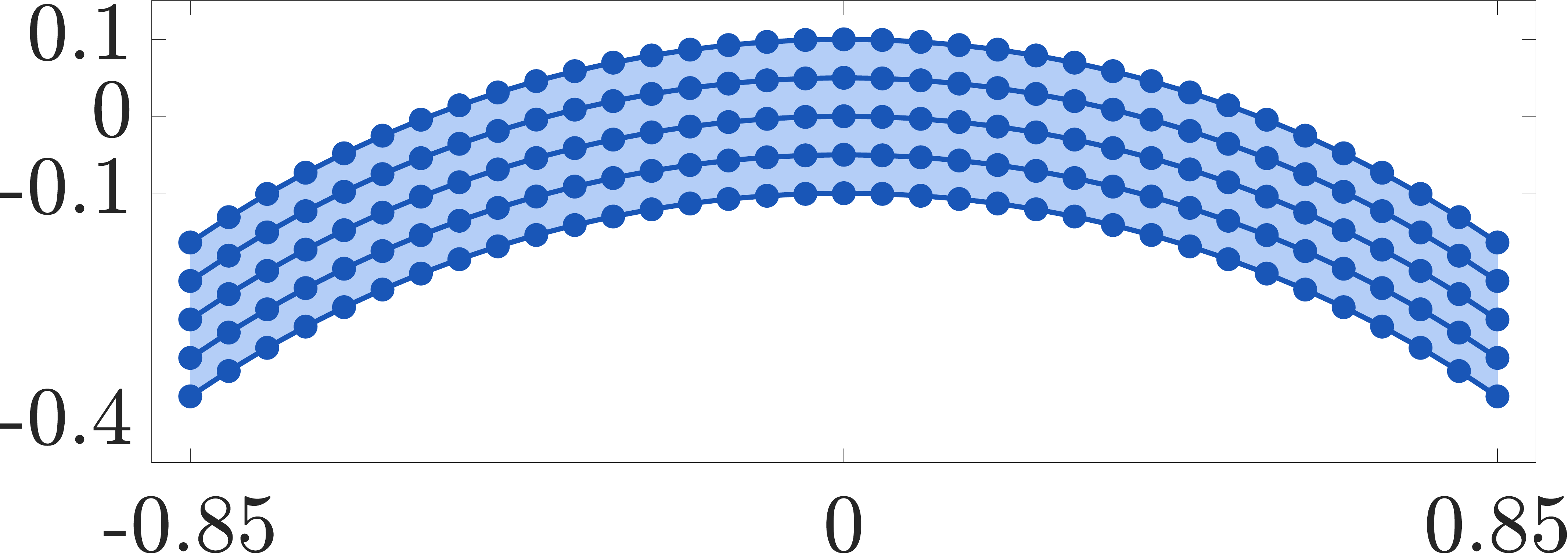}
			\caption{Undeformed shape}
			\label{fig:cap_template}
		\end{subfigure}
		\vspace{10pt}
		
		\begin{subfigure}{\textwidth}
			\centering
			\includegraphics[width = \textwidth]{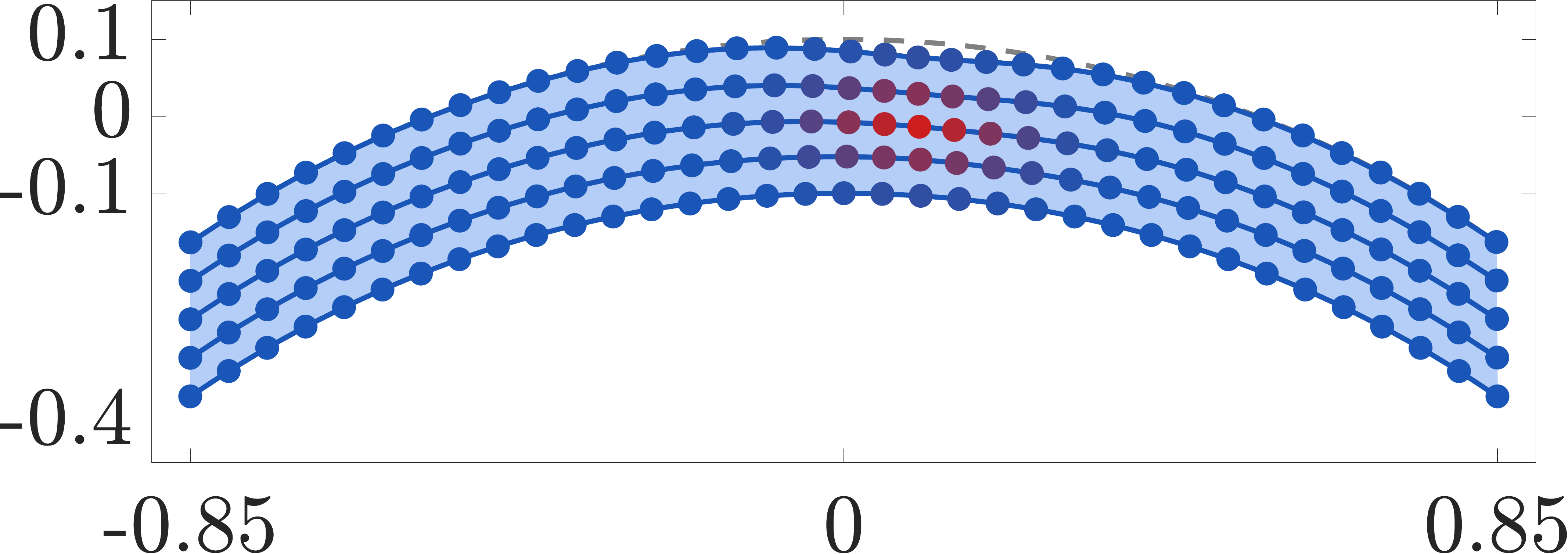}
			\caption{Simulated target}
			\label{fig:cap_target}
		\end{subfigure}
	\end{subfigure}
	\hfill
	\begin{subfigure}{0.53\textwidth}
		\centering
		\includegraphics[width = \textwidth]{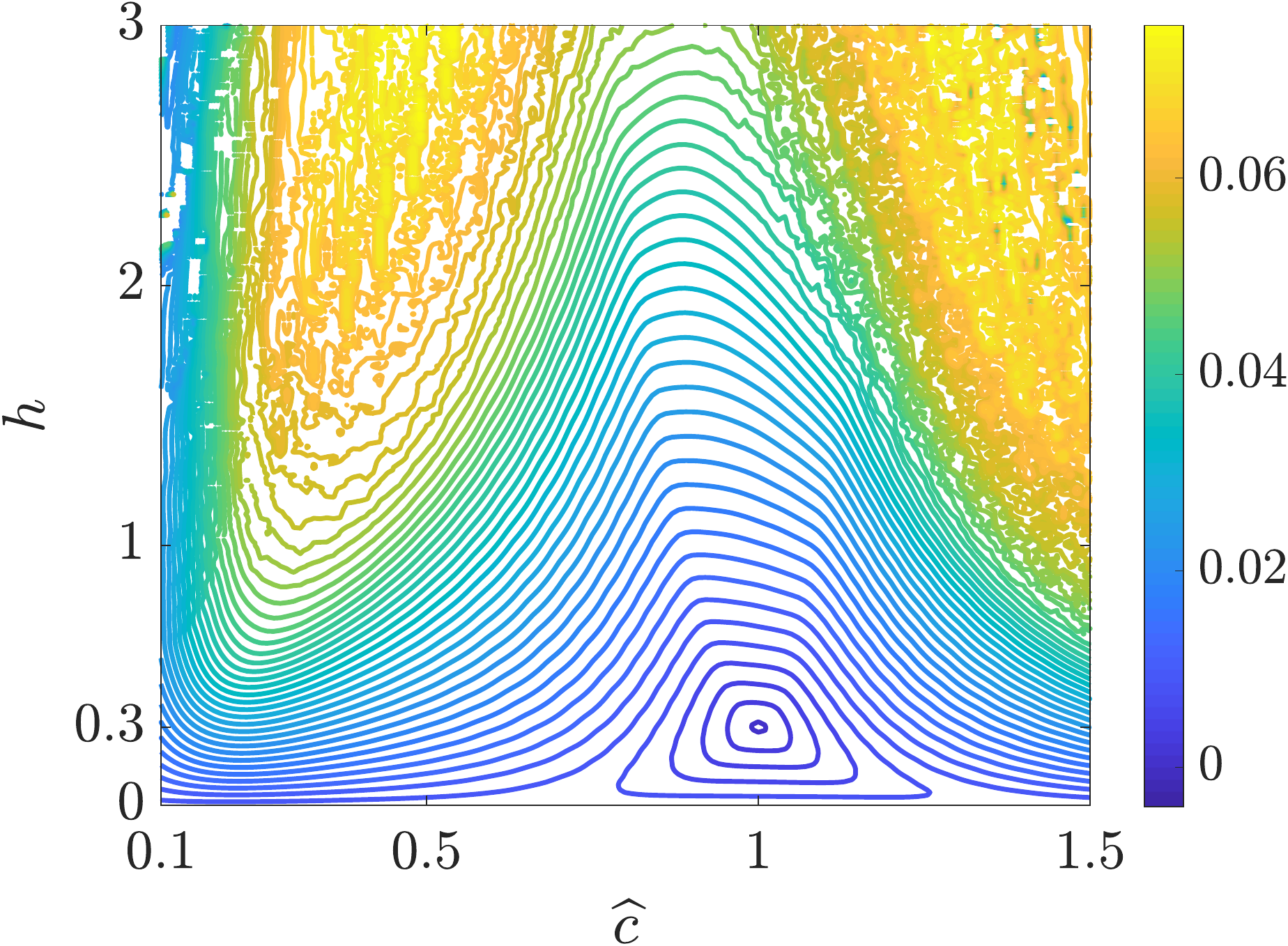}
		\caption{Objective function}
		\label{fig:cap_obj}
	\end{subfigure}
	\caption{2D simulated experiment with the true solution at $(\widehat c, h) = (1, 0.3)$}
	\label{fig:cap_inverse}
\end{figure}

\begin{figure}
	\centering
	\begin{subfigure}{0.43\textwidth}
		\begin{subfigure}{\textwidth}
			\centering
			\includegraphics[width = \textwidth]{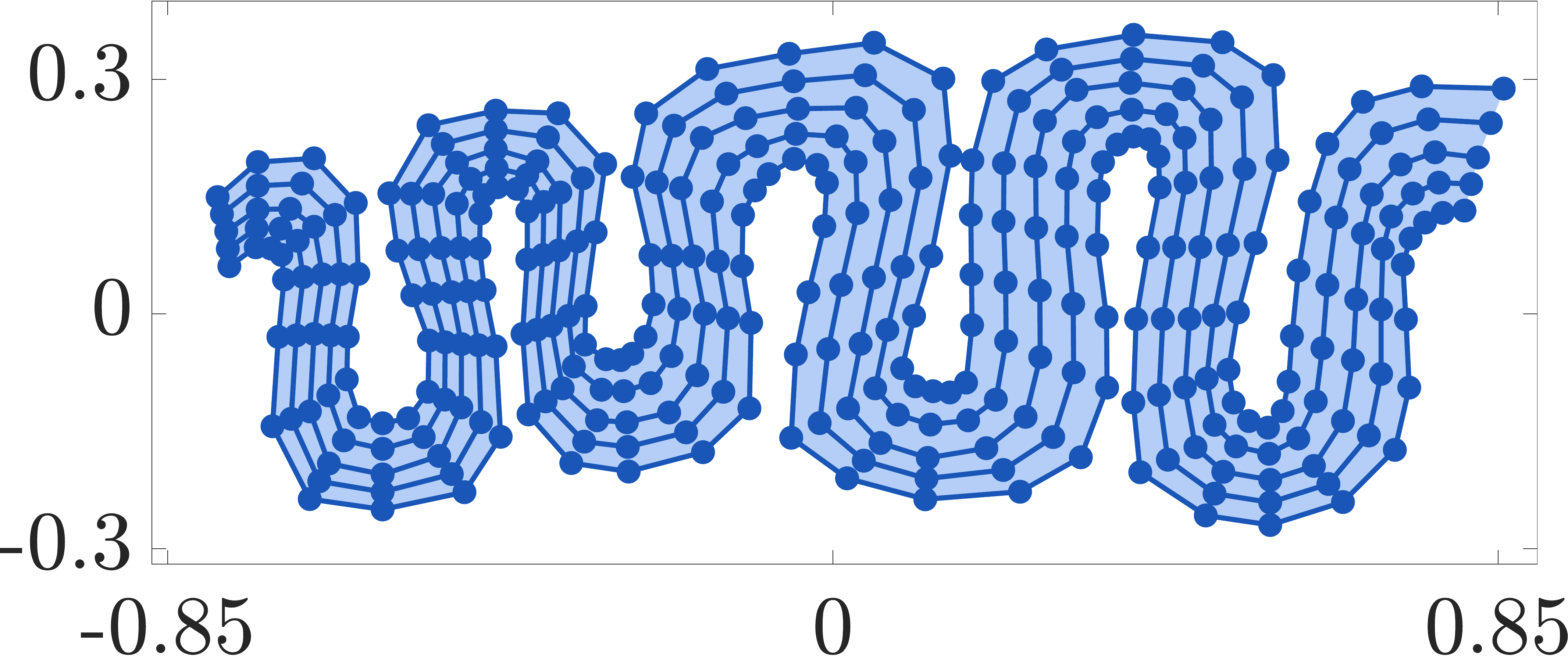}
			\caption{Undeformed shape}
			\label{fig:fold_template}
		\end{subfigure}
		\vspace{10pt}
		
		\begin{subfigure}{\textwidth}
			\centering
			\includegraphics[width = \textwidth]{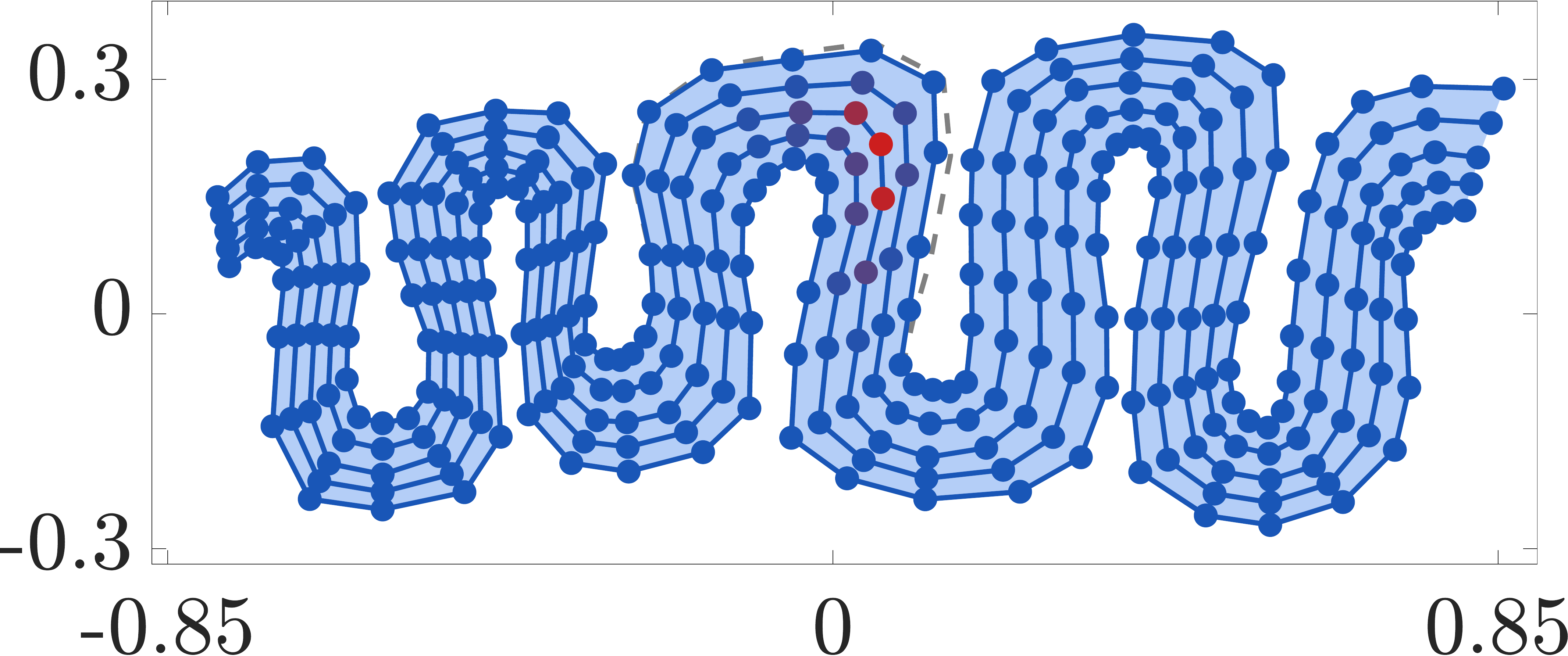}
			\caption{Simulated target}
			\label{fig:fold_target}
		\end{subfigure}
	\end{subfigure}
	\hfill
	\begin{subfigure}{0.53\textwidth}
		\centering
		\includegraphics[width = \textwidth]{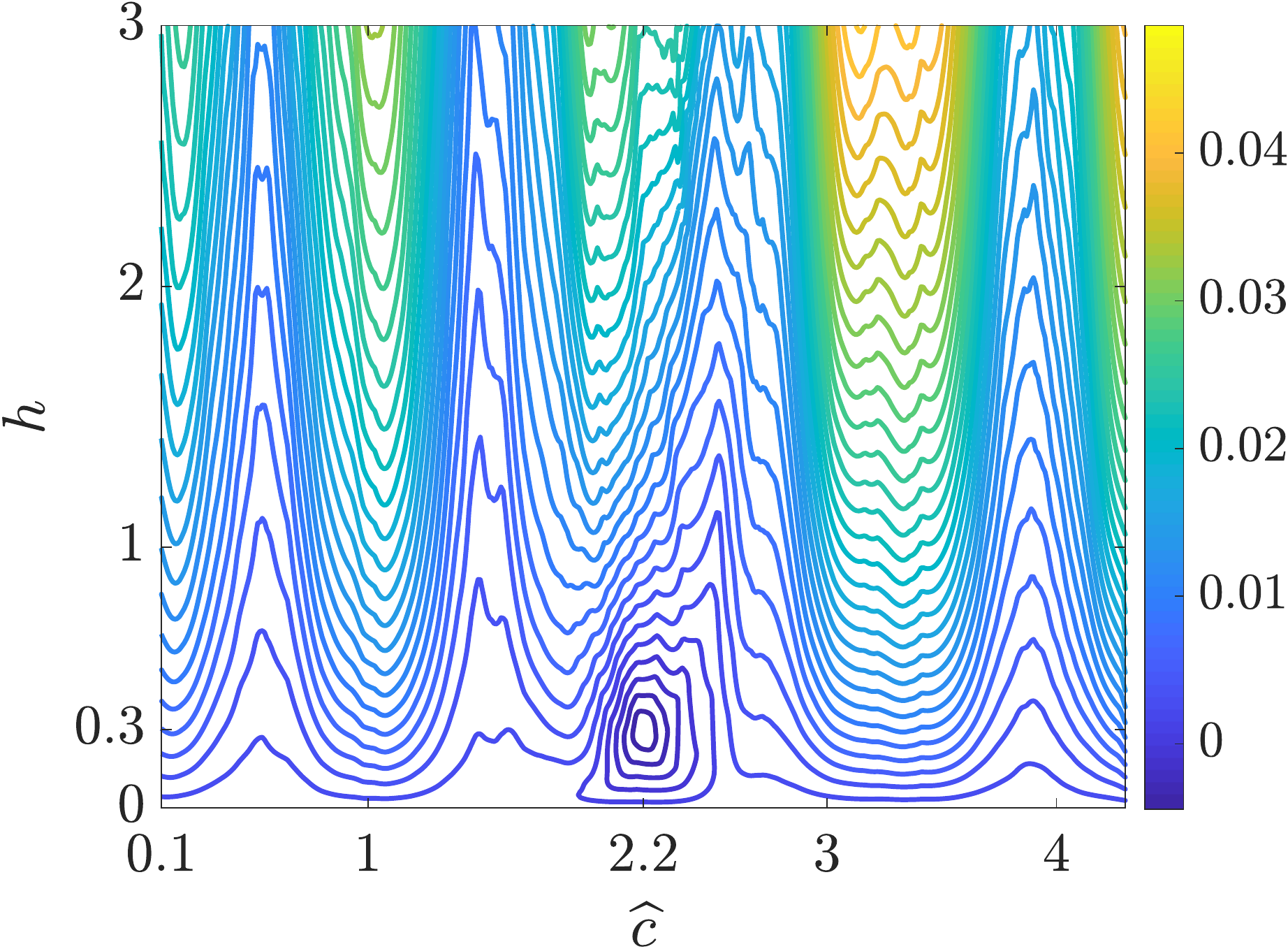}
		\caption{Objective function}
		\label{fig:fold_obj}
	\end{subfigure}
	\caption{2D simulated experiment with the true solution at $(\widehat c, h) = (2.2, 0.3)$}
	\label{fig:fold_inverse}
\end{figure}

Two-dimensional results are shown in Figs.~\ref{fig:cap_inverse} and~\ref{fig:fold_inverse}. We used a regular shape in Fig.~\ref{fig:cap_inverse}, while Fig.~\ref{fig:fold_inverse} was created to resemble a lateral cut of cerebral cortex. The initial shapes are given in Figs.~\ref{fig:cap_template} and ~\ref{fig:fold_template}, and the simulated targets are shown in Figs.~\ref{fig:cap_target} and ~\ref{fig:fold_target}, in which the dashed curves indicate the positions of the initial shapes.
We used a very small force in both cases for the purpose of demonstration. Recall that $c$ is restricted on the middle layer, and therefore only has one degree of freedom, which we use the arc length from the leftmost landmark in both cases as its coordinate representation $\widehat c$. The level curves of the objective function $\widehat J(\widehat c, h)$ are presented in Figs.~\ref{fig:cap_obj}~and~\ref{fig:fold_obj}, which demonstrate that the global minimizers are very close to the force parameters generating the targets, even though the deformation is small. Fig.~\ref{fig:fold_obj} also suggests that the optimization problem could be challenging for a highly curved shape. For example, if a gradient descent-like algorithm starts at $(\widehat c, h) = (1, 1)$, it would be difficult to locate the minimizer around $(2.2, 0.3)$.

\begin{wrapfigure}[8]{r}{0.5\textwidth}
	\centering
	\includegraphics[width = 0.45\textwidth]{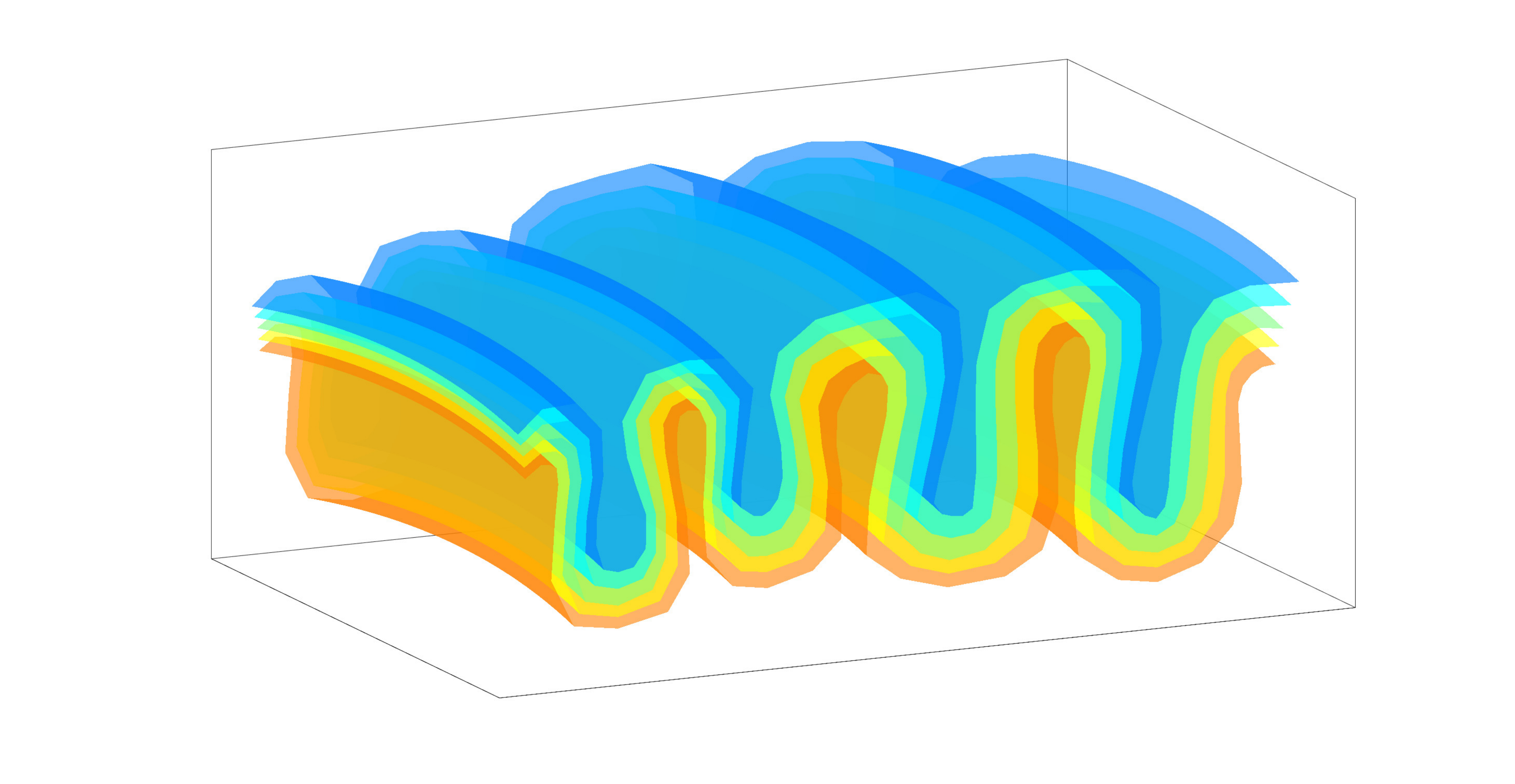}
	\caption{A 3D synthetic folded shape.}
	\label{fig:3D} 
\end{wrapfigure}

Our 3D experiments use the initial shape shown in Fig.~\ref{fig:3D}, which is a synthetic cerebral cortex.
Table~\ref{tbl:3d_results} compares the estimated parameters $(\widehat c, h)$ with the ground truth for a series of simulated deformations. Again, we see that the true solution for a simulated target can be accurately retrieved by our method.

\begin{table}
	\setlength{\tabcolsep}{2pt}
	\caption{Estimated force parameters compared to their ground truth.}
	\label{tbl:3d_results}
	\vspace{3pt}
	\centering
	\begin{tabular}{c cccc}
		\toprule \\[-13.8pt]
		\rowcolor{Gray}
		\vphantom{\rule[-5pt]{10pt}{15pt}}
		$(\widehat c_{\mathrm{true}}, h_{\mathrm{true}})$ &
		(1.60, 0.50, 0.30) &
		(2.20, 0.50, 0.30) &
		(1.60, 0.65, 0.30) &
		(2.20, 0.65, 0.30) \\[-2pt]
		\cmidrule{1-5}
		$(\widehat c, h)$ &
		(1.60, 0.50, 0.31) &
		(2.21, 0.49, 0.29) &
		(1.60, 0.65, 0.32) &
		(2.21, 0.65, 0.29) \\
		\midrule \\[-13.8pt]
		\rowcolor{Gray}
		\vphantom{\rule[-5pt]{10pt}{15pt}}
		$(\widehat c_{\mathrm{true}}, h_{\mathrm{true}})$ &
		(1.60, 0.50, 0.35) &
		(2.20, 0.50, 0.35) &
		(1.60, 0.65, 0.35) &
		(2.20, 0.65, 0.35) \\[-2pt]
		\cmidrule{1-5}
		$(\widehat c, h)$ &
		(1.60, 0.49, 0.36) &
		(2.20, 0.50, 0.34) &
		(1.60, 0.65, 0.36) &
		(2.21, 0.65, 0.34) \\
		\bottomrule
	\end{tabular}
\end{table}


\section{Discussion}
\label{sec:discussion}
We have presented in this paper a new model describing shape evolution where the control can be interpreted as the derivative of a body force density in the deforming volume. We also have provided a preliminary set of experiments, using simulated data, and based on derivative-free optimization methods. 

Current and future work include applying this approach to medical imaging data, and in particular to the study of atrophy due to neurodegenerative diseases. This will require a more general definition of the force field than the single sink model we have used, and an interesting problem will be to maintain parametric identifiability for such models, by injecting suitable biological priors. Using more sophisticated optimization methods will require being able to explicitly compute gradients, which, even though formally feasible, represent serious numerical challenges, probably involving a solution of the linear system $L_\varphi v = j\,dp$ on parallel hardware.

We have not discussed in this paper any theoretical result on the well-posedness of the considered problems, such as sufficient conditions on the existence of solutions of \eqref{eq:syst.1}, or consistency of the discretization schemes. We however included what we believe to be suitable smoothness assumptions regarding, in particular, the RKHS $V$. Such results are under investigation and will be published in the near future.  


\bibliographystyle{splncs04}
\bibliography{Elasticity}

\end{document}